\input amssym.def
\input amssym.tex

\headline={\ifnum \pageno=1 {\hfill} 
\else{\hss \tenrm -- \folio\ -- \hss}\fi}
\footline={\hfil}

\def\dater{\vglue-10mm\rightline{(\the\day/\the\month/\the\year)}}

\hsize 146mm
\vsize 224mm
\hoffset=6mm
\voffset=8mm
\baselineskip=5mm
\overfullrule =0pt

 at 10,5pt

\font\GGgtitre=cmbx10 at 17pt

%Debut symboles speciaux

\def\og{\leavevmode\raise.30ex
\hbox{$\scriptscriptstyle\langle\!\langle\>$}}    
\def\bigog{\leavevmode\raise.30ex
\hbox{$\langle\!\langle\>$}}  
%guillemets ouvrants
\def\fg{\leavevmode\raise.30ex
\hbox{$\scriptscriptstyle\>\rangle\!\rangle$}}    
\def\bigfg{\leavevmode\raise.30ex
\hbox{$\>\rangle\!\rangle$}}    
%guillemets fermants

\catcode`\@=11

\font\author=cmcsc10
\font\pauthor=cmcsc10 at 8pt
\font\tenmsx=msam10
\font\sevenmsx=msam10 scaled 700
\font\fivemsx=msam10 scaled 500
\font\tenmsy=msbm10
\font\sevenmsy=msbm10 scaled 700
\font\fivemsy=msbm10 scaled 500
\newfam\msxfam
\newfam\msyfam
\textfont\msxfam=\tenmsx  \scriptfont\msxfam=\sevenmsx
\scriptscriptfont\msxfam=\fivemsx
\textfont\msyfam=\tenmsy  \scriptfont\msyfam=\sevenmsy
\scriptscriptfont\msyfam=\fivemsy

\def\hexnumber@#1{\ifnum#1<10 \number#1\else
\ifnum#1=10 A\else\ifnum#1=11 B\else\ifnum#1=12 C\else
\ifnum#1=13 D\else\ifnum#1=14 E\else\ifnum#1=15 F\fi\fi\fi\fi\fi\fi\fi}

\def\msx@{\hexnumber@\msxfam}
\def\msy@{\hexnumber@\msyfam}
\mathchardef\nmid="3\msy@2D
\mathchardef\varnothing="0\msy@3F
\mathchardef\nexists="0\msy@40
\mathchardef\smallsetminus="2\msy@72
\def\Bbb{\ifmmode\let\next\Bbb@\else
\def\next{\errmessage{Use \string\Bbb\space only in math mode}}\fi\next}
\def\Bbb@#1{{\Bbb@@{#1}}}
\def\Bbb@@#1{\fam\msyfam#1}

\font\tentbl=cmr10 scaled 900
\font\seventbl=cmr7 scaled 900
\font\fivetbl=cmr5 scaled 900

\newfam\tblfam

\textfont\tblfam=\tentbl
\scriptfont\tblfam=\seventbl
\scriptscriptfont\tblfam=\fivetbl

% Fin symboles speciaux

\def \Z {{\Bbb Z}}

\font\ci= eufm10

\def\ciU{\hbox{\ci U}}

 at 9,5pt

\font\f=cmr10 at 7pt

\def\f1{\hbox{\f 1}}
\def\f2{\hbox{\f 2}}
\def\f3{\hbox{\f 3}}
\def\f4{\hbox{\f 4}}
\def\f5{\hbox{\f 5}}
\def\f6{\hbox{\f 6}}
\def\f7{\hbox{\f 7}}
\def\f8{\hbox{\f 8}}
\def\f9{\hbox{\f 9}}

\def \d {\,{\rm d}}

\def\le{\leqslant}
\def\ge{\geqslant}

\topskip=10pt
\font\sept=cmti9

\def\rightheadline{\ifnum\pageno=\chstart{\hfill}
          \else{\centerline{\sept Chen's double sieve, 
Goldbach's conjecture and the twin prime problem, 2}}
          \hfill \hskip -6mm \tenrm\folio\fi}
\def\leftheadline{\ifnum\pageno=\chstart{\hfill}
          \else\tenrm\folio \hskip -3,5mm \hfill{\centerline{\pauthor 
J. Wu}}\fi}
\headline={\ifnum\pageno=\chstart{\hfill}
\else{\ifodd\pageno\rightheadline\else\leftheadline\fi}\fi}
\footline={\hfill}

\pageno=1
\newcount\chstart
\chstart=\pageno

%\dater

\vglue 5mm

\centerline{\GGgtitre Chen's double sieve, Goldbach's conjecture}

\bigskip

\centerline{\GGgtitre and the twin prime problem, 2}

\vskip 8mm

\centerline{\author J. Wu}

\vskip 8mm

{\leftskip=12mm
\rightskip=12mm
{{\bf Abstract}.
For every even integer $N$,
denote by $D_{1,2}(N)$ the number of representations of
$N$ as a sum of a prime and an integer having at most two prime
factors.
In this paper, we give a new lower bound for $D_{1,2}(N)$.
\par}}

\footnote{{}}{2000 {\it Mathematics Subject Classification}: 
11P32, 11N35, 11N05.}

\vskip 10mm

{\hsize 122mm
\baselineskip=4mm
\overfullrule =0pt
\leftskip=22mm

\par}

\vskip 5mm

\noindent{\bf \S\ 1. Introduction}

\medskip

Let $\Omega(n)$ be the number of all prime factors of the integer $n$ 
with the convention $\Omega(1)=0$.
For each even integer $N\ge 4$,
we define 
$$D(N) := |\{p\le N :\Omega(N - p) = 1\}|,$$
where and in what follows,
the letter $p$, with or without subscript, denotes a prime number.
The well known Goldbach conjecture can be stated as $D(N)\ge 1$ for
all even integers $N\ge 4$.
A more precise version of this conjecture was proposed by Hardy \&
Littlewood [10]:
$$D(N) \sim 2 \Theta(N)
\qquad(N\to \infty),
\leqno(1.1)$$
where
$$C_N:=\prod_{p\mid N, \, p>2} {p-1\over p-2} \prod_{p>2}
\bigg(1-{1\over (p-1)^2}\bigg),
\qquad
\Theta(N) := {C_N N\over (\log N)^2}.
\leqno(1.2)$$
Certainly, the asymptotic formula (1.1) is extremely difficult.
One way of approaching the lower bound problem in (1.1) is
to give a non-trivial lower bound for the quantity
$$D_{1,2}(N) := |\{p\le N : \Omega(N-p)\le 2\}|.$$
In this direction, Chen [5] proved, by his system of
weights and the switching principle,
the following famous theorem:
{\sl Every sufficiently large even integer can be written as sum of a
prime and an integer
having at most two prime factors}.
More precisely he established
$$D_{1,2}(N)\ge 0.67 \, \Theta(N)
\leqno(1.3)$$
for $N\ge N_0$.
As Halberstam \& Richert indicated in [9], 
it would be interesting to know whether a
more elaborate weighting procedure
could be adapted to the purpose of (1.3).
This might lead to numerical improvements and could be important.
Chen's constant 0.67 has been improved by many authors.
The historical record is as follows:
$$\eqalign{
& \hbox{$0.689$ by Halberstam \& Richert [9]},
\cr
& \hbox{$0.754$ by Chen [6]},
\cr
& \hbox{$0.81$ \hskip 1,75mm by Chen [7]},
\cr
& \hbox{$0.828$ by Cai \& Lu [4]},
\cr
& \hbox{$0.836$ by Wu [13]},
\cr
& \hbox{$0.867$ by Cai [2]}.
\cr}$$

The aim of this paper is to propose a better constant.

\proclaim Theorem.
For sufficiently large $N$, we have
$$D_{1,2}(N)\ge 0.899\,\Theta(N).$$

Our improvement comes from a delicate application of Chen's double sieve 
([8], [12], [13]),
which can be described as follows:
With standard notation in theory of sieve method,
the linear sieve formulas (see [9], or Lemma 2.2 of [13])
can be stated as
$$X V(z) f\bigg({\log Q\over \log z}\bigg) + {\rm error}
\le S({\cal A}; {\cal P}, z)
\le X V(z) F\bigg({\log Q\over \log z}\bigg) + {\rm error}.
\leqno(1.4)$$
These inequalities are the best possible in the sense that taking
$${\cal A}
= {\cal B}_\nu
:= \{n\le x : \Omega(n)\equiv \nu\,({\rm mod}\,2)\}
\qquad(\nu=1, 2),$$
the upper and  lower bounds in (1.4) are respectively attained
by $\nu=1$ and $\nu=2$ (see [9], page 239).
Aiming at a better Bombieri-Davenport's upper bound [1]
$$D(N)\le \{8+o(1)\}\Theta(N),$$
Chen [8] found improvement for  (1.4) 
for some special sequences ${\cal A}$.
Roughly speaking, for the sequence
$${\cal A} = \{N-p : p\le N\}$$
he narrowed down the gap in (1.4) 
by introducing two functions $h(s)$ and $H(s)$ 
such that the functions $sf(s)/(2e^\gamma)$ and $sF(s)/(2e^\gamma)$ 
are replaced by $sf(s)/(2e^\gamma)+h(s)$ and $sF(s)/(2e^\gamma)-H(s)$  
respectively,
where $\gamma$ is the Euler constant.
The key point is thus to prove $h(s)>0$ and $H(s)>0$.
Chen's proof is very long and somewhat difficult
to follow, but his innovative idea is clear (see [11] for example).
In [13], we gave a more comprehensive treatment on this method
and name it as Chen's double sieve. Indeed, 
our treatment is not only simpler but even more powerful than Chen's.
Our approach improved Chen's upper estimate
$D(N)\le 7.8342\Theta(N)$ to $D(N)\le 7.8209\Theta(N)$.
It is worth to indicate that Chen's record stood for 26 years 
before our work [13].

To prove our Theorem, we first simplify and 
improve Chen's weight system (compare (12) of [7] and Lemma 2.2 below), 
and then apply Chen's double sieve,
as the classical linear sieve, 
to handle terms such as 
$\Upsilon_2$,
$\Upsilon_3$, 
$\Upsilon_4$, 
$\Upsilon_5$ and $\Upsilon_6$
(cf. Propositions 4.1, 4.2, 4.3 and 4.4 below). 
The idea of using Chen's double sieve to treat sums of the type
$$\sum_{\scriptstyle N^{\phi_1}\le 
p<N^{\phi_2}\atop\scriptstyle (p, N)=1}
S({\cal A}_p; {\cal P}(N), N^\kappa)
\leqno(1.5)$$
was first appeared in [12]. 
However, due to the first condition in (3.1) below, 
a direct application of our Chen's double sieve 
can only handle the initial part of the sum over small $p$ in (1.5)
(i.e. $p\le N^{1/4}$). On the other hand, 
very recently Cai [2] used a similar idea to 
control the sum over large $p$ in (1.5). 
Actually his method can be viewed as a simplified version 
of Chen's double sieve 
(see Proposition 4.4 below and the comments before it).
Here we shall combine both versions  
and refine them to obtain our result.
Apparently from the proof, 
we shall see that the first version gives a saving of 0.0211 
while the second saves 0.0078.
Without Chen's double sieve technique, 
we still obtain 0.870 in place of 0.899,
which is slightly better than Cai's 0.867.

Clearly our method can be used to refine the corresponding constants
in the conjugate problems ([2] and [3]).
The proofs are very similar and even easier and simpler.
Hence we omit the relevant discussion.
Maybe this is a good exercise for senior graduate students 
in analytic number theory.

\vskip 5mm

\noindent{\bf \S\ 2. Chen's system of weights}

\medskip

This section is devoted to discuss the weighted sieve of Chen type.
Let 
$${\cal A} := \{N-p : p\le N\}
\qquad{\rm and}\qquad
{\cal P}(N) := \{p : (p, N) = 1\}.$$
The sieve function is defined as
$$
S({\cal A}; {\cal P}(N), z)
:= |\{a\in {\cal A} : (a, P(z))=1\}|,$$
where $P(z):=\prod_{p\le z, \, p\in {\cal P}(N)} p$.

\proclaim Lemma 2.1.
Let $0<\kappa<\sigma\le {1\over 3}$. Then we have
$$2D_{1,2}(N)\ge
2S({\cal A}; {\cal P}(N), N^\kappa)
- S_1(\kappa, \sigma)
- 2S_2(\kappa, \sigma)
- S_3(\kappa, \sigma)
+ S_4(\kappa, \sigma)
+ O(N^{1-\kappa}),
\leqno(2.1)$$
where 
$$\eqalign{
S_1(\kappa, \sigma)
& :=\sum_{\scriptstyle N^\kappa\le p<N^\sigma\atop\scriptstyle (p, N)=1}
S({\cal A}_p; {\cal P}(N), N^\kappa),
\cr
S_2(\kappa, \sigma)
& := \mathop{\sum \,\, \sum}_{\scriptstyle N^\sigma\le 
p_1<p_2<(N/p_1)^{1/2}
\atop\scriptstyle (p_1p_2, N)=1}
S({\cal A}_{p_1 p_2}; {\cal P}(Np_1), p_2),
\cr
S_3(\kappa, \sigma) 
& := \mathop{\sum \,\, \sum}_{\scriptstyle N^\kappa\le 
p_1<N^\sigma\le p_2<(N/p_1)^{1/2}
\atop\scriptstyle (p_1p_2, N)=1}
S({\cal A}_{p_1 p_2}; {\cal P}(Np_1), p_2),
\cr
S_4(\kappa, \sigma)
& := \mathop{\sum \,\, \sum \,\, \sum}_{\scriptstyle N^\kappa\le 
p_1<p_2<p_3<N^\sigma
\atop\scriptstyle (p_1p_2p_3, N)=1}
S({\cal A}_{p_1 p_2 p_3}; {\cal P}(Np_1), p_2).
\cr}$$

The inequality (2.1) first appeared in [7] (page 479, (11)) with 
$(\kappa, \sigma) 
= ({1\over 12}, {1\over 3.047}), ({1\over 9.2}, {1\over 3.41})$
without proof. 
Cai \& [Lu] [4] gave a proof with an extra assumption 
$3\sigma+\kappa>1$.
In [13], we proved (2.1) 
under the hypothesis $0<\kappa<\sigma<{1\over 3}$.
Clearly the proof there is also valid for $\sigma={1\over 3}$.
Very recently Cai [2] gave another proof for Lemma 2.1.

As in [7], we shall apply (2.1) with 
two different pairs of parameters
$(\kappa, \sigma)$ to take advantage of $S_4(\kappa, \sigma)$.
Our weighted sieve is simpler and more poweful than those of Chen 
([7], (12)) and Cai ([2], Lemma 6).

\proclaim Lemma 2.2.
Let $\kappa_2>\kappa_1\ge {1/18}$
such that 
$3\kappa_1+\kappa_2<1/2$
and 
$3\kappa_1-\kappa_2<1/6$.
Then we have
$$4D_{1,2}(N)
\ge
3\Upsilon_1
+ \Upsilon_2
- \Upsilon_3
- \Upsilon_4
+ \Upsilon_5
+ \Upsilon_6
- 2\Upsilon_7
- \Upsilon_8
- \Upsilon_9
- \Upsilon_{10}
- \Upsilon_{11}
+ O(N^{1-\kappa_1}),
\leqno(2.2)$$
where
$$\eqalign{
& \Upsilon_i 
:= S({\cal A}; {\cal P}(N), N^{\kappa_i})
\quad(i=1, 2),
\cr\noalign{\vskip 5mm}
& \Upsilon_3 
:= \sum_{\scriptstyle N^{\kappa_1}\le 
p<N^{1/3}\atop\scriptstyle (p, N)=1}
S({\cal A}_p; {\cal P}(N), N^{\kappa_1}),
\cr
& \Upsilon_4 
:= \sum_{\scriptstyle N^{\kappa_1}\le p<N^{1/2-3\kappa_1}
\atop\scriptstyle (p, N)=1}
S({\cal A}_p; {\cal P}(N), N^{\kappa_1}),
\cr
& \Upsilon_5 := \mathop{\sum \,\, \sum }_{\scriptstyle 
N^{\kappa_1}\le p_1<p_2<N^{\kappa_2}
\atop\scriptstyle (p_1p_2,N)=1}
S({\cal A}_{p_1p_2}; {\cal P}(N), N^{\kappa_1}),
\cr
& \Upsilon_6 := \mathop{\sum \,\, \sum}_{\scriptstyle
N^{\kappa_1}\le p_1<N^{\kappa_2}\le p_2<N^{1/2-3\kappa_1}
\atop\scriptstyle (p_1p_2,N)=1}
S({\cal A}_{p_1p_2}; {\cal P}(N), N^{\kappa_1}),
\cr
& \Upsilon_7 := \mathop{\sum \,\, \sum}_{
\scriptstyle N^{1/2-3\kappa_1}\le p_1<p_2<(N/p_1)^{1/2}
\atop{\atop\scriptstyle (p_1p_2, N)=1}}
S({\cal A}_{p_1 p_2}; {\cal P}(Np_1), p_2),
\cr
& \Upsilon_8 :=
\mathop{\sum \,\, \sum}_{\scriptstyle N^{\kappa_1}\le 
p_1<N^{1/3}\le p_2<(N/p_1)^{1/2}
\atop\scriptstyle (p_1p_2, N)=1}
S({\cal A}_{p_1 p_2}; {\cal P}(Np_1), p_2),
\cr
& \Upsilon_9
:= \mathop{\sum \,\, \sum}_{\scriptstyle 
N^{\kappa_2}\le p_1<N^{1/2-3\kappa_1}\le p_2<(N/p_1)^{1/2}
\atop\scriptstyle (p_1p_2, N)=1}
S\big({\cal A}_{p_1 p_2}; {\cal P}(Np_1), (N/p_1p_2)^{1/2}\big),
\cr
& \Upsilon_{10}
:= \mathop{\sum \,\, \sum \,\, \sum \,\, \sum}_{\scriptstyle
N^{\kappa_1}\le p_1<p_2<p_3<p_4<N^{\kappa_2} \atop\scriptstyle 
(p_1p_2p_3p_4,N)=1}
S({\cal A}_{p_1p_2p_3p_4}; {\cal P}(N), p_2),
\cr
& \Upsilon_{11}
:= \mathop{\sum \,\, \sum \,\, \sum \,\, \sum}_{\scriptstyle
N^{\kappa_1}\le p_1<p_2<p_3<N^{\kappa_2}\le p_4<N^{1/2-2\kappa_1}/p_3
\atop\scriptstyle (p_1p_2p_3p_4,N)=1}
S({\cal A}_{p_1p_2p_3p_4}; {\cal P}(N), p_2).
\cr}$$

\noindent{\sl Proof}.
By noticing that our hypothesis
implies $\kappa_2<1/2-3\kappa_1\le 1/3$,
we can apply (2.1) with $(\kappa, \sigma)=(\kappa_2, 1/2-3\kappa_1)$ 
to obtain
$$2 D_{1,2}(N)
\ge
2 \Upsilon_2 
- S_1(\kappa_2, 1/2-3\kappa_1)
- 2 \Upsilon_7
- S_3(\kappa_2, 1/2-3\kappa_1)
+ O(N^{1-\kappa_2}),
\leqno(2.3)$$
where the term $S_4(\kappa_2, 1/2-3\kappa_1)$ is dropped by non-negativity.

Buchstab's identity, when applied three times, gives the equality
$$\eqalign{\Upsilon_2
& = \Upsilon_1
- \sum_{\scriptstyle N^{\kappa_1}\le p<N^{\kappa_2}
\atop\scriptstyle (p,N)=1}
S({\cal A}_p; {\cal P}(N), N^{\kappa_1})
+ \Upsilon_5
- \mathop{\sum \,\, \sum \,\, \sum}_{\scriptstyle N^{\kappa_1}\le 
p_1<p_2<p_3<N^{\kappa_2}
\atop\scriptstyle (p_1p_2p_3,N)=1}
S({\cal A}_{p_1p_2p_3}; {\cal P}(N), p_1).
\cr}$$
Similarly,
a twice application of Buchstab's identity yields
$$\eqalign{
S_1(\kappa_2,  1/2-3\kappa_1)
& = \sum_{\scriptstyle N^{\kappa_2}\le p<N^{1/2-3\kappa_1}
\atop\scriptstyle (p,N)=1}
S({\cal A}_p; {\cal P}(N), N^{\kappa_1})
- \Upsilon_6
\cr
& \quad
+ \mathop{\sum \,\, \sum \,\, \sum}_{\scriptstyle 
N^{\kappa_1}\le p_1<p_2<N^{\kappa_2}\le p_3<N^{1/2-3\kappa_1}
\atop\scriptstyle (p_1p_2p_3,N)=1}
S({\cal A}_{p_1p_2p_3}; {\cal P}(N), p_1).
\cr}$$
By Buchstab's identity, we can prove
$$\eqalign{S_3(\kappa_2,  1/2-3\kappa_1)
& \le \Upsilon_9
+ \mathop{\sum \,\, \sum \,\, \sum}_{\scriptstyle
N^{\kappa_2}\le p_1<N^{1/2-3\kappa_1}\le 
p_2<p_3<(N/p_1p_2)^{1/2}\atop\scriptstyle (p_1p_2p_3, N)=1}
S({\cal A}_{p_1 p_2 p_3}; {\cal P}(Np_1), p_3).
\cr}$$
Inserting them into (2.3), we find that
$$2 D_{1,2}(N)
\ge \Upsilon_1
+ \Upsilon_2
- \Upsilon_4
+ \Upsilon_5
+ \Upsilon_6
- 2 \Upsilon_7
- \Upsilon_9
- \Delta_1
+ O(N^{1-\kappa_2}),
\leqno(2.4)$$
where
$$\eqalign{
\Delta_1
& := \mathop{\sum \,\, \sum \,\, \sum}_{\scriptstyle 
N^{\kappa_1}\le p_1<p_2<p_3<N^{\kappa_2}
\atop\scriptstyle (p_1p_2p_3,N)=1}
S({\cal A}_{p_1p_2p_3}; {\cal P}(N), p_1) 
\cr
& \quad
+ \mathop{\sum \,\, \sum \,\, \sum}_{\scriptstyle N^{\kappa_1}\le 
p_1<p_2<N^{\kappa_2}\le p_3<N^{1/2-3\kappa_1}
\atop\scriptstyle (p_1p_2p_3,N)=1}
S({\cal A}_{p_1p_2p_3}; {\cal P}(N), p_1)
\cr
& \quad
+ \mathop{\sum \,\, \sum \,\, \sum}_{\scriptstyle
N^{\kappa_2}\le p_1<N^{1/2-3\kappa_1}\le 
p_2<p_3<(N/p_1p_2)^{1/2}\atop\scriptstyle (p_1p_2p_3, N)=1}
S({\cal A}_{p_1 p_2 p_3}; {\cal P}(Np_1), p_3).
\cr}$$

The inequality (2.1) with $(\kappa, \sigma) = (\kappa_1, {1/3})$ 
gives 
$$2 D_{1,2}(N)
\ge
2 \Upsilon_1 
- \Upsilon_3
- \Upsilon_8
+ S_4(\kappa_1, {1/3})
+ O(N^{1-\kappa_1}),
\leqno(2.5)$$
where we have used the fact that $S_2(\kappa_1, {1/3})=0$.

Adding (2.4) to (2.5) yields
$$4 D_{1,2}(N)
\ge
3 \Upsilon_1
+ \Upsilon_2 
- \Upsilon_3
- \Upsilon_4
+ \Upsilon_5
+ \Upsilon_6
- 2 \Upsilon_7
- \Upsilon_8
- \Upsilon_9
+ \Delta_2
+ O(N^{1-\kappa_1}),
\leqno(2.6)$$
where
$$
\Delta_2
:= \mathop{\sum \,\, \sum \,\, \sum}_{\scriptstyle N^{\kappa_1}\le 
p_1<p_2<p_3<N^{1/3}
\atop\scriptstyle (p_1p_2p_3, N)=1}
S({\cal A}_{p_1 p_2 p_3}; {\cal P}(N), p_2)
- \Delta_1.
$$

Clearly all the summation ranges in the three triple sums of  
$\Delta_1$ are distinct
and the first two are covered 
in the range of the triple sum
in $\Delta_2$ 
(since our hypothesis on $\kappa_1$ and $\kappa_2$ implies
$\max\{\kappa_2, 1/2-3\kappa_1\}\le 1/3$).
On the other hand, we easily see that the range of summation
in the third triple sum of $\Delta_1$ is equivalent to
$N^{\kappa_2}\le p_1<N^{1/2-3\kappa_1}\le p_2\le (N/p_1)^{1/3}$ 
and $p_2<p_3<(N/p_1p_2)^{1/2}$.
From this we deduce that
$(N/p_1p_2)^{1/2}
\le N^{(1/2+3\kappa_1-\kappa_2)/2}
\le N^{1/3}$,
since $3\kappa_1-\kappa_2<1/6$.
Thus this range is also contained 
in the triple sum of $\Delta_2$. 
Therefore we have
$$\eqalign{\Delta_2
& \ge - \mathop{\sum \,\, \sum \,\, \sum}_{\scriptstyle 
N^{\kappa_1}\le p_1<p_2<p_3<N^{\kappa_2}
\atop\scriptstyle (p_1p_2p_3,N)=1}
\big\{S({\cal A}_{p_1p_2p_3}; {\cal P}(N), p_1) - S({\cal 
A}_{p_1p_2p_3}; {\cal P}(N), p_2)\big\}
\cr
& \hskip 4,8mm
- \mathop{\sum \,\, \sum \,\, \sum}_{\scriptstyle
N^{\kappa_1}\le p_1<p_2<N^{\kappa_2}\le p_3<N^{1/2-2\kappa_1}/p_2
\atop\scriptstyle (p_1p_2p_3,N)=1}
\big\{S({\cal A}_{p_1p_2p_3}; {\cal P}(N), p_1) - S({\cal 
A}_{p_1p_2p_3}; {\cal P}(N), p_2)\big\}
\cr
& \hskip 4,8mm
+ \mathop{\sum \,\, \sum \,\, \sum}_{\scriptstyle
N^{\kappa_2}\le p_1<N^{1/2-3\kappa_1}\le p_2<p_3<(N/p_1p_2)^{1/2}
\atop\scriptstyle (p_1p_2p_3, N)=1}
\big\{S({\cal A}_{p_1p_2p_3}; {\cal P}(N), p_2) - S({\cal 
A}_{p_1p_2p_3}; {\cal P}(N), p_3)\big\}
\cr
& \ge - \Upsilon_{10} - \Upsilon_{11} + O(N^{1-\kappa_1}).
\cr}$$
Combining with (2.6), we obtain the required result.
\hfill
$\square$

\smallskip

{\bf Remark 1}.
Apparently from the proof,
we have choosen 
$(\kappa, \sigma)
=(\kappa_1, 1/2-3\kappa_1),\,
(\kappa_2, 1/3)$ in the application of Lemma 2.1.
It is possible to optimize the choice of $\sigma$.
But this augments the number of terms of (2.2) 
and the numeric improvement for Theorem is quite small.

\vskip 5mm

\noindent{\bf \S\ 3. Chen's double sieve}

\medskip

In this section, 
we recall Chen's double sieve described in [13]
and give numeric lower bounds for $H(s)$ and $h(s)$ for later use.

For any large even integer $N$, we write
$${\cal A} := \{N - p : p\le N\},
\qquad
{\cal P}(N) := \{p : (p, N) = 1\}.$$
Let $\delta>0$ be a sufficiently small number
\footnote{$^{(*)}$}
{In numerical computation, we can formally take $\delta=0$.} 
and $k\in \Z$.
Put 
$$Q := N^{1/2-\delta},
\qquad
\underline d := Q/d,
\qquad
{\cal L} := \log N,
\qquad
W_k := N^{\delta^{1+k}}.$$
Denote by $\pi_{[Y, Z)}$ the characteristic function of the set
${\cal P}(N) \cap [Y, Z)$.
For $k\in \Z^+$ and $N\ge 2$,
let $\ciU_k(N)$ be the set of all arithmetical functions $\sigma$ which 
can be written as the form
$$\sigma = \pi_{[V_1/\Delta, V_1)}*\cdots*\pi_{[V_i/\Delta, V_i)},$$
where $\Delta$ is a real number with
$1 + {\cal L}^{-4}\le \Delta<1 + 2 {\cal L}^{-4}$,
$i$ is an integer with $0\le i\le k$, and $V_1, \dots, V_i$ are 
real numbers satisfying
$$\cases{
V_1^2\le Q,                            & {}
\cr\noalign{\smallskip}
V_1 V_2^2\le Q,                        & {}
\cr\noalign{\smallskip}
\cdots\cdots\cdots\cdots\cdot          & {}
\cr\noalign{\smallskip}
V_1\cdots V_{i-1} V_i^2\le Q,          & {}
\cr\noalign{\smallskip}
V_1\ge V_2\ge \cdots\ge V_i\ge W_k.    & {}
\cr}
\leqno(3.1)$$
We adopt the convention that $\sigma$ is the characteristic function of the set 
$\{1\}$ if $i=0$.

Let $F$ and $f$ be defined by 
$$\eqalign{
F(s) = 2e^\gamma/s,
&
\quad\qquad
f(s)=0
\qquad
(0<s\le 2),
\cr
(sF(s))' = f(s-1),
&
\qquad\quad
(sf(s))' = F(s-1)
\quad
(s>2),
\cr}
\leqno(3.2)$$
where $\gamma$ is Euler's constant.
Moreover we take
$$A(s) := sF(s)/2e^\gamma,
\qquad
a(s) := sf(s)/2e^\gamma,
\leqno(3.3)$$
and introduce the notation
$$\leqalignno{
\Phi(N, \sigma, s) 
& := \sum_d \sigma(d) S({\cal A}_d; {\cal P}(d N), 
\underline d^{1/s}),
& (3.4)
\cr
\Theta(N, \sigma)
& := 4 {\rm li}(N)
\sum_d {\sigma(d) C_{d N}\over \varphi(d)\log \underline d},
& (3.5)
\cr}$$
where $\varphi(d)$ is the Euler function.

For $k\in \Z^+$, $N_0\ge 2$ and $s\in [1, 10]$, we define
$H_{k,N_0}(s)$ and $h_{k,N_0}(s)$ to be the supremum of $h\ge -\infty$
such that for all $N\ge N_0$ and $\sigma\in \ciU_k(N)$, 
the  inequalities
$$\eqalign{
\Phi(N, \sigma, s)
& \le \{A(s) - h\} \, \Theta(N, \sigma),
\cr\noalign{\vskip 1mm}
\Phi(N, \sigma, s)
& \ge \{a(s) + h\} \, \Theta(N, \sigma)
\cr}$$
hold true respectively.
Obviously $H_{k,N_0}(s)$ and $h_{k,N_0}(s)$ are decreasing in $N_0$,
as well as decreasing in $k$ by Lemma 3.1.
Hence their limits at infinity exist (in the extended real line), and we write
$$\eqalign{
H_k(s)
& := \lim_{N_0\rightarrow \infty} H_{k,N_0}(s),
\cr
H(s)
& := \lim_{k\rightarrow \infty} H_k(s),
\cr}
\qquad
\eqalign{
h_k(s)
& := \lim_{N_0\rightarrow \infty} h_{k,N_0}(s),
\cr
h(s)
& := \lim_{k\rightarrow \infty} h_k(s).
\cr}$$

The next lemma collects the concerned properties of these functions
(see [13], Lemma 3.2, Propositions 1 \& 2 and Corollary 1).

\proclaim Lemma 3.1.
(i)
For $k\in \Z^+, N\ge N_0, s\in [1, 10]$ and $\sigma\in \ciU_k(N)$,
we have
$$\leqalignno{\Phi(N, \sigma, s)
& \le \{A(s) - H_{k,N_0}(s)\} \Theta(N, \sigma),
& (3.6)
\cr\noalign{\vskip 1mm}
\Phi(N, \sigma, s)
& \ge \{a(s) + h_{k,N_0}(s)\} \Theta(N, \sigma).
& (3.7)
\cr}$$

\vskip -1,8mm
{\sl 
(ii)
For $k\in \Z^+$ and $s\in [1,10]$, we have 
$H_k(s)\ge 0$ and $h_k(s)\ge 0$.

(iii)
For $2\le s\le s'\le 10$, we have
$$h(s)\ge h(s') + \int_{s-1}^{s'-1} {H(t)\over t} \d t
\qquad\hbox{and}\qquad
H(s)\ge H(s') + \int_{s-1}^{s'-1} {h(t)\over t} \d t.
\leqno(3.8)$$

(iv)
The function $H(s)$ is decreasing on $[1, 10]$.
The function $h(s)$ is increasing on $[1, 2]$ 
and is decreasing on $[2, 10]$.}

\smallskip

We cannot give explicit expressions for $H(s)$ and $h(s)$.
But it is tractable to obtain numeric lower bounds 
for these two functions.
Let
$$s_i:=2+0.1\times i\quad(i\ge 0).
\leqno(3.9)$$
By ([13], \S\ 7),
we have the numeric lower bounds of $H(s_i)$ for $2\le i\le 10$.
Next we shall consider the case of $11\le i\le 29$ 
and the lower bounds of $h(s_i)$ for $0\le i\le 29$.
These will be used in the proof of Theorem.

Let ${\bf 1}_{[a, b]}(t)$ be the characteristic function of the 
interval $[a, b]$ and
$$\sigma(a, b, c)
:= \int_a^b \log\bigg({c\over t-1}\bigg) {\d t\over t},
\qquad
\sigma_0(t) 
:= {\sigma(3, t+2, t+1)\over 1-\sigma(3, 5, 4)}.$$
From (6.2) of [13] and the decreasing property of $H(s)$, 
we deduce 
$$H(s_j)\ge \sum_{2\le i\le 10} c_{i, j} H(s_i),
\leqno(3.10)$$
for $11\le j\le 29$,
where
$$\eqalign{
c_{2, j}
& := \int_{1}^{s_2}
\bigg\{{\sigma_0(t)\over t}\log\bigg({4\over s_j-1}\bigg)
+{{\bf 1}_{[s_j-2,3]}(t)\over t}\log\bigg({t+1\over s_j-1}\bigg)\bigg\}
\d t,
\cr
c_{i, j}
& := \int_{s_{i-1}}^{s_i}
\bigg\{{\sigma_0(t)\over t}\log\bigg({4\over s_j-1}\bigg)
+{{\bf 1}_{[s_j-2,3]}(t)\over t}\log\bigg({t+1\over s_j-1}\bigg)\bigg\}
\d t
\quad
(3\le i\le 10).
\cr}$$
From the first inequality of (3.8) and the fact that $h(s)\ge 0$, 
we also derive
$$\leqalignno{h(s_j)
& \ge \int_{s_j-1}^5 {H(t)\over t} \d t
& (3.11)
\cr
& \ge H(s_{2})\log\bigg({s_{\max\{2, j-10\}}\over s_j-1}\bigg) 
+\sum_{\max\{3, j-9\}\le i\le 29} 
H(s_i)\log\bigg({s_i\over s_{i-1}}\bigg)
\cr}$$
for $0\le j\le 29$.

Using the numeric lower bounds of $H(s_i)$ for $2\le i\le 10$
given in ([13], \S\ 7), (3.10) and (3.11), we get via a
numerical computation the following results.
$$\displaylines{
\vbox{\tabskip = 0pt\offinterlineskip
\halign{
\vrule # & &\hfil$ $ $#$ $ \!\! $ \hfil & \vrule #\cr
\noalign{\hrule}
height 2mm
&&&&&&&&&&&&&&&&&&
\cr
& \,\,     i       \,\,               &
& \,\,\,   s_i     \,\,\,             &
& \hskip 2,5mm H(s_i)\ge \hskip 2,5mm &
& \,\,\,     i     \,\,\,             &
& \,\,\,   s_i     \,\,\,             &
& \hskip 2,5mm H(s_i)\ge \hskip 2,5mm &
& \,\,\,     i     \,\,\,             &
& \,\,\,   s_i     \,\,\,             &
& \hskip 2,5mm H(s_i)\ge \hskip 2,5mm &
\cr
height 2mm
&&&&&&&&&&&&&&&&&&
\cr
\noalign{\hrule}
height 2mm
&&&&&&&&&&&&&&&&&&
\cr
&           &
&           &
&           &
& 10        &
& 3.0       &
& 0.0072943 &
& 20        &
& 4.0       &
& 0.0010835 &
\cr
height 2mm
&&&&&&&&&&&&&&&&&&
\cr
\noalign{\hrule}
height 2mm
&&&&&&&&&&&&&&&&&&
\cr
&           &
&           &
&           &
& 11        &
& 3.1       &
& 0.0061642 &
& 21        &
& 4.1       &
& 0.0008451 &
\cr
height 2mm
&&&&&&&&&&&&&&&&&&
\cr
\noalign{\hrule}
height 2mm
&&&&&&&&&&&&&&&&&&
\cr
& 2         &
& 2.2       &
& 0.0223939 &
& 12        &
& 3.2       &
& 0.0052233 &
& 22        &
& 4.2       &
& 0.0006482 &
\cr
height 2mm
&&&&&&&&&&&&&&&&&&
\cr
\noalign{\hrule}
height 2mm
&&&&&&&&&&&&&&&&&&
\cr
& 3         &
& 2.3       &
& 0.0217196 &
& 13        &
& 3.3       &
& 0.0044073 &
& 23        &
& 4.3       &
& 0.0004882 &
\cr
height 2mm
&&&&&&&&&&&&&&&&&&
\cr
\noalign{\hrule}
height 2mm
&&&&&&&&&&&&&&&&&&
\cr
& 4         &
& 2.4       &
& 0.0202876 &
& 14        &
& 3.4       &
& 0.0036995 &
& 24        &
& 4.4       &
& 0.0003602 &
\cr
height 2mm
&&&&&&&&&&&&&&&&&&
\cr
\noalign{\hrule}
height 2mm
&&&&&&&&&&&&&&&&&&
\cr
& 5         &
& 2.5       &
& 0.0181433 &
& 15        &
& 3.5       &
& 0.0030860 &
& 25        &
& 4.5       &
& 0.0002592 &
\cr
height 2mm
&&&&&&&&&&&&&&&&&&
\cr
\noalign{\hrule}
height 2mm
&&&&&&&&&&&&&&&&&&
\cr
& 6         &
& 2.6       &
& 0.0158644 &
& 16        &
& 3.6       &
& 0.0025551 &
& 26        &
& 4.6       &
& 0.0001803 &
\cr
height 2mm
&&&&&&&&&&&&&&&&&&
\cr
\noalign{\hrule}
height 2mm
&&&&&&&&&&&&&&&&&&
\cr
& 7         &
& 2.7       &
& 0.0129923 &
& 17        &
& 3.7       &
& 0.0020972 &
& 27        &
& 4.7       &
& 0.0001187 &
\cr
height 2mm
&&&&&&&&&&&&&&&&&&
\cr
\noalign{\hrule}
height 2mm
&&&&&&&&&&&&&&&&&&
\cr
& 8         &
& 2.8       &
& 0.0100686 &
& 18        &
& 3.8       &
& 0.0017038 &
& 28        &
& 4.8       &
& 0.0000702 &
\cr
height 2mm
&&&&&&&&&&&&&&&&&&
\cr
\noalign{\hrule}
height 2mm
&&&&&&&&&&&&&&&&&&
\cr
& 9         &
& 2.9       &
& 0.0078162 &
& 19        &
& 3.9       &
& 0.0013680 &
& 29        &
& 4.9       &
& 0.0000313 &
\cr
height 2mm
&&&&&&&&&&&&&&&&&&
\cr
\noalign {\hrule}}}
\cr}
$$
\vskip -3mm
\centerline{Table 1. Numeric lower bounds for $H(s_i)$}

\medskip

$$\displaylines{
\vbox{\tabskip = 0pt\offinterlineskip
\halign{
\vrule # & &\hfil$ $ $#$ $ \!\! $ \hfil & \vrule #\cr
\noalign{\hrule}
height 2mm
&&&&&&&&&&&&&&&&&&
\cr
& \,\,     i       \,\,           &
& \,\,\,   s_i     \,\,\,         &
& \hskip 3mm h(s_i)\ge \hskip 3mm &
& \,\,\,     i     \,\,\,         &
& \,\,\,   s_i     \,\,\,         &
& \hskip 3mm h(s_i)\ge \hskip 3mm &
& \,\,\,     i     \,\,\,         &
& \,\,\,   s_i     \,\,\,         &
& \hskip 3mm h(s_i)\ge \hskip 3mm &
\cr
height 2mm
&&&&&&&&&&&&&&&&&&
\cr
\noalign{\hrule}
height 2mm
&&&&&&&&&&&&&&&&&&
\cr
& 0         &
& 2.0       &
& 0.0232385 &
& 10        &
& 3.0       &
& 0.0077162 &
& 20        &
& 4.0       &
& 0.0010120 &
\cr
height 2mm
&&&&&&&&&&&&&&&&&&
\cr
\noalign{\hrule}
height 2mm
&&&&&&&&&&&&&&&&&&
\cr
& 1         &
& 2.1       &
& 0.0211041 &
& 11        &
& 3.1       &
& 0.0066236 &
& 21        &
& 4.1       &
& 0.0008099 &
\cr
height 2mm
&&&&&&&&&&&&&&&&&&
\cr
\noalign{\hrule}
height 2mm
&&&&&&&&&&&&&&&&&&
\cr
& 2         &
& 2.2       &
& 0.0191556 &
& 12        &
& 3.2       &
& 0.0055818 &
& 22        &
& 4.2       &
& 0.0006440 &
\cr
height 2mm
&&&&&&&&&&&&&&&&&&
\cr
\noalign{\hrule}
height 2mm
&&&&&&&&&&&&&&&&&&
\cr
& 3         &
& 2.3       &
& 0.0173631 &
& 13        &
& 3.3       &
& 0.0046164 &
& 23        &
& 4.3       &
& 0.0005084 &
\cr
height 2mm
&&&&&&&&&&&&&&&&&&
\cr
\noalign{\hrule}
height 2mm
&&&&&&&&&&&&&&&&&&
\cr
& 4         &
& 2.4       &
& 0.0157035 &
& 14        &
& 3.4       &
& 0.0037529 &
& 24        &
& 4.4       &
& 0.0003980 &
\cr
height 2mm
&&&&&&&&&&&&&&&&&&
\cr
\noalign{\hrule}
height 2mm
&&&&&&&&&&&&&&&&&&
\cr
& 5         &
& 2.5       &
& 0.0141585 &
& 15        &
& 3.5       &
& 0.0030123 &
& 25        &
& 4.5       &
& 0.0003085 &
\cr
height 2mm
&&&&&&&&&&&&&&&&&&
\cr
\noalign{\hrule}
height 2mm
&&&&&&&&&&&&&&&&&&
\cr
& 6         &
& 2.6       &
& 0.0127132 &
& 16        &
& 3.6       &
& 0.0023901 &
& 26        &
& 4.6       &
& 0.0002365 &
\cr
height 2mm
&&&&&&&&&&&&&&&&&&
\cr
\noalign{\hrule}
height 2mm
&&&&&&&&&&&&&&&&&&
\cr
& 7         &
& 2.7       &
& 0.0113556 &
& 17        &
& 3.7       &
& 0.0018997 &
& 27        &
& 4.7       &
& 0.0001791 &
\cr
height 2mm
&&&&&&&&&&&&&&&&&&
\cr
\noalign{\hrule}
height 2mm
&&&&&&&&&&&&&&&&&&
\cr
& 8         &
& 2.8       &
& 0.0100756 &
& 18        &
& 3.8       &
& 0.0015336 &
& 28        &
& 4.8       &
& 0.0001336 &
\cr
height 2mm
&&&&&&&&&&&&&&&&&&
\cr
\noalign{\hrule}
height 2mm
&&&&&&&&&&&&&&&&&&
\cr
& 9         &
& 2.9       &
& 0.0088648 &
& 19        &
& 3.9       &
& 0.0012593 &
& 29        &
& 4.9       &
& 0.0000981 &
\cr
height 2mm
&&&&&&&&&&&&&&&&&&
\cr
\noalign {\hrule}}}
\cr}
$$
\vskip -3mm
\centerline{Table 2. Numeric lower bounds for $h(s_i)$}

\bigskip

{\bf Remark 2}.
It is possible to get better numeric lower bounds for $H(s_i)$ and
$h(s_i)$ by applying (3.8) repeatedly.
But the improvement will be small.

\vskip 5mm

\noindent{\bf \S\ 4. Application of Chen's double sieve}

\medskip

In this section, 
we apply Chen's double sieve to estimate 
the terms $\Upsilon_3$, $\Upsilon_4$,
$\Upsilon_5$ and $\Upsilon_6$ in (2.2).
Propositions 4.1, 4.2, 4.3 and 4.4 below are results in general context.
These estimates are better than those obtained 
by the classical linear sieve,
since $H(s)>0$ and $h(s)>0$.

\proclaim Proposition 4.1.
Let $0<\phi_1<\phi_2<1/4$ and $\kappa>0$ 
such that $\phi_2+\kappa\le 1/2$.
Then for $N\to\infty$, we have
$$\sum_{\scriptstyle N^{\phi_1}\le p<N^{\phi_2}
\atop\scriptstyle (p, N)=1}
S({\cal A}_{p}; {\cal P}(N), N^\kappa)
\le \bigg\{8
\int_{(1/2-\phi_2)/\kappa}^{(1/2-\phi_1)/\kappa} 
{A(t)-H(t)\over t(1-2\kappa t)}\d t
+o(1)\bigg\}\Theta(N).
$$

\noindent{\sl Proof}.
We keep use of the previous notation.
Denote by $S$ the sum in the proposition.
Let 
$\alpha_j:=N^{\phi_1}\Delta^j$
and $J$ be the integer such that
$\alpha_J\le N^{\phi_2}<\alpha_{J+1}$.
We write 
$$S
= \sum_{1\le j\le J} \sum_{p}
\pi_{[\alpha_{j-1}, \alpha_j)}(p)
S\big({\cal A}_{p}; {\cal P}(pN), \underline{p}^{1/\tau_p}\big)
+ R_1,
\leqno(4.1)$$
where $\tau_p := (\log\underline{p})/(\kappa\log N)$ 
and
$$R_1
:= \sum_{\alpha_J\le p<N^{\phi_2}} 
S({\cal A}_{p}; {\cal P}(N), N^\kappa)
\ll \sum_{\alpha_J\le p<N^{\phi_2}} {N/p}
\ll \Theta(N){\cal L}^{-3}.
\leqno(4.2)$$

Introducing
$$\tau_j:=(\log\underline{\alpha_j})/(\kappa\log N),$$
we easily see that
$\pi_{[\alpha_{j-1}, \alpha_j)}(p)
\not=0
\Rightarrow
\tau_j\le \tau_p\le \tau_{j-1}$.
Thus we can deduce from (4.1) and (4.2) that
$$S
\le
\sum_{1\le j\le J} \sum_{p}
\pi_{[\alpha_{j-1}, \alpha_j)}(p)
S\big({\cal A}_{p}; {\cal P}(pN), \underline{p}^{1/\tau_j}\big)
+ O\big(\Theta(N){\cal L}^{-3}\big),
\leqno(4.3)$$
where we have used the following estimates:
$$\eqalign{ 
& \sum_{1\le j\le J} \sum_{p}
\pi_{[\alpha_{j-1}, \alpha_j)}(p)
\big\{
S\big({\cal A}_{p}; {\cal P}(pN), \underline{p}^{1/\tau_p}\big)
- S\big({\cal A}_{p}; {\cal P}(pN), \underline{p}^{1/\tau_j}\big)
\big\}
\cr
& \quad
\le \sum_{1\le j\le J} \sum_{\alpha_{j-1}\le p<\alpha_j}
\sum_{\underline{p}^{1/\tau_p}\le p'<\underline{p}^{1/\tau_j}} 
N/(pp')
\cr
& \quad
\ll N{\cal L}^{-5}
\sum_{1\le j\le J} \sum_{\alpha_{j-1}\le p<\alpha_j} 1/p
\cr
& \quad
\ll \Theta(N){\cal L}^{-3}.
\cr}$$

Next we treat the inner sum (over $p$) in (4.3).
Clearly for each $j\in \{1, \dots, J\}$,
our hypothesis on $\phi_1, \phi_2$ and $\kappa$
assures that
the function $\pi_{[\alpha_{j-1}, \alpha_j)}\in \ciU_k(N)$
for all $k\ge 0$, $N_0\ge 2$ and $N\ge N_0$,
and $\tau_j\ge 1$.
Thus we can apply (3.6) of Lemma 3.1 to estimate the sum over $p$
(which is $\Phi(N, \pi_{[\alpha_{j-1}, \alpha_j)}, \tau_j)$) :
$$\eqalign{S
& \le \sum_{1\le j\le J}
\{A(\tau_j) - H_{k, N_0}(\tau_j)\}
\Theta(N, \pi_{[\alpha_{j-1}, \alpha_j)}) 
+ O\big(\Theta(N){\cal L}^{-3}\big)
\cr
& \le 4 {\rm li}(N)
{C_N\over \log\underline{1}}
\sum_{\alpha_0\le p<\alpha_J}
{A(\tau_p)-H_{k, N_0}(\tau_p)\over (p-2)(1-\log p/\log \underline{1})}
+ O\big(\Theta(N){\cal L}^{-3}\big)
\cr
& \le 4 {\rm li}(N)
{C_N\over \log\underline{1}}
\sum_{N^{\phi_1}\le p<N^{\phi_2}}
{A(\tau_p) - H_{k, N_0}(\tau_p)\over
(p - 2) (1 - \log p/\log \underline{1})}
+ O\big(\Theta(N){\cal L}^{-3}\big),
\cr}$$
where we have used the fact that $A(s)-H_{k, N_0}(s)$ 
is increasing in $s$.
An integration by parts with the prime number theorem shows that
$$\sum_{N^{\phi_1}\le p<N^{\phi_2}}
{A(\tau_p) - H_{k, N_0}(\tau_p)\over
(p - 2) (1 - \log p/\log \underline{1})}
= \int_{(1/2-\phi_2)/\kappa}^{(1/2-\phi_1)/\kappa}
{A(t) - H_{k, N_0}(t)\over t(1-2\kappa t)}\d t 
+ O_{\delta, k}(\varepsilon).$$
Hence
$$S
\le 8\bigg\{
\int_{(1/2-\phi_2)/\kappa}^{(1/2-\phi_1)/\kappa} 
{A(t) - H_{k, N_0}(t)\over t(1-2\kappa t)}\d t
+ O_{\delta, k}(\varepsilon)\bigg\}
\Theta(N)
$$
for $N\ge N_0$.
From this, we infer that
$$\limsup_{N\to\infty} {S\over \Theta(N)}
\le 8\int_{(1/2-\phi_2)/\kappa}^{(1/2-\phi_1)/\kappa} 
{A(t) - H_{k, N_0}(t)\over t(1-2\kappa t)}\d t
+ O_{\delta, k}(\varepsilon),
$$
which implies, by taking $N\to\infty$, $k\to\infty$ 
and $\varepsilon\to 0$,
$$
\limsup_{N\to\infty} {S\over \Theta(N)}
\le 8\int_{(1/2-\phi_2)/\kappa}^{(1/2-\phi_1)/\kappa} 
{A(t) - H(t)\over t(1-2\kappa t)}\d t.
$$
Clearly this is equivalent to the required inequality.
\hfill 
$\square$

\goodbreak
\medskip

In a similar fashion we can prove the following results.

\proclaim Proposition 4.2.
Let $0<\phi_1<\phi_2<1/6$ and $\kappa>0$ 
such that $2\phi_2+\kappa\le 1/2$.
Then for $N\to\infty$, we have
$$\eqalign{
& \mathop{\sum\;\sum}_{\scriptstyle N^{\phi_1}\le p_1<p_2<N^{\phi_2}
\atop\scriptstyle (p_1p_2, N)=1}
S({\cal A}_{p_1p_2}; {\cal P}(N), N^\kappa)
\cr
& \hskip 6mm
\ge \bigg\{8 
\int_{\phi_1}^{\phi_2} 
\int_{(1/2-\phi_2-t)/\kappa}^{(1/2-2t)/\kappa}
{a(u)+h(u)\over t(1-2t-2\kappa u)} \d t\d u
+o(1)\bigg\}\Theta(N).
\cr}$$

\proclaim Proposition 4.3.
Let $0<\phi_1<\phi_2\le \phi_3<\phi_4<1/4$ and $\kappa>0$ 
such that $2\phi_2+\phi_4<1/2$ 
and $\phi_2+\phi_4+\kappa\le 1/2$.
Then for $N\to\infty$, we have
$$\eqalign{
& \mathop{\sum_{N^{\phi_1}\le p_1<N^{\phi_2}}\;
\sum_{N^{\phi_3}\le p_1<N^{\phi_4}}}_{(p_1p_2, N)=1}
S({\cal A}_{p_1p_2}; {\cal P}(N), N^\kappa)
\cr
& \hskip 10mm
\ge \bigg\{8
\int_{\phi_1}^{\phi_2} 
\int_{(1/2-\phi_4-t)/\kappa}^{(1/2-\phi_3-t)/\kappa}
{a(u)+h(u)\over tu(1-2t-2\kappa u)}\d t\d u
+o(1)\bigg\}\Theta(N).
\cr}$$

\smallskip

Finally we estimate the sum of the type in (1.5) with $\phi_1\ge 1/4$.
In this case, we cannot directly apply our delicate Chen's double sieve 
because of the first condition of (3.1).
As what Cai [2] remarked, 
it is possible to use a simplified version of Chen's double sieve.
This approach will give a result 
better than using the classic linear sieve 
but weaker than Proposition 4.1,
since, without iteration, 
$\Psi_1(s)$ or $\Psi_2(s)$ are principal contributions of $H(s)$.
(See Lemmas 5.1 and 5.2 of [13]
and compare Proposition 4.4 below and Proposition 4.1.)

\proclaim Proposition 4.4.
Let $\kappa>0$, $\phi>0$ and $2\le s\le 3\le s'\le 5$ 
such that $1/4\le 1/2-s\kappa<\phi$.
Then for $N\to\infty$, we have
$$\sum_{\scriptstyle N^{1/2-s\kappa}\le p<N^{\phi}
\atop\scriptstyle (p, N)=1}
S({\cal A}_{p}; {\cal P}(N), N^{\kappa})
\le \bigg\{8
\int_{(1/2-\phi)/\kappa}^s {A(t)-\Psi_1(s)\over t(1-2\kappa t)}\d t
+ o(1)\bigg\}
\Theta(N),
$$
where
$$\eqalign{
\Psi_1(s)
& := -\int_2^{s'-1} {\log(t-1)\over t}\d t
+ {1\over 2} \int_{1-1/s}^{1-1/s'} 
{\log(s't-1)\over t(1-t)} \d t
\cr
& \hskip 4,5mm
- \max_{\phi\ge 2} 
\mathop{\int\int\int}_{1/s'\le t\le u\le v\le 1/s}
\omega\bigg({\phi - t - u - v\over u}\bigg)
{\d t \d u \d v\over t u^2 v}
\cr}$$
and $\omega(u)$ is Buchstab's function.
The same result also holds if we replace $\Psi_1(s)$ by $\Psi_2(s)$,
where the function $\Psi_2(s)$ is defined as in Lemma 5.2 of [13]. 

\noindent{\sl Proof}.
For simplicity, we denote the sum by $S$.
Since $N^\kappa\ge \underline{p}^{1/s}$ for $p\ge N^{1/2-s\kappa}$,
we can write
$$\eqalign{S
& \le \sum_{\scriptstyle N^{1/2-s\kappa}\le p<N^{\phi}
\atop\scriptstyle (p, N)=1}
S\big({\cal A}_{p}; {\cal P}(N), \underline{p}^{1/s}\big)
\cr
& \le \sum_{1\le j\le J} 
\sum_p\pi_{[\alpha_{j-1}, \alpha_j)}(p)
S\big({\cal A}_{p}; {\cal P}(N), \underline{p}^{1/s}\big),
\cr}$$
where $\alpha_j:=N^{1/2-s\kappa}\Delta^j$ and 
$J$ is the integer such that $\alpha_{J-1}\le N^\phi<\alpha_J$.

Similar to Lemma 4.1 of [13], 
we can prove that there is a constant $\eta>0$ such that
$$S\le \sum_{1\le j\le J} 
\sum_p\pi_{[\alpha_{j-1}, \alpha_j)}(p)
\bigg(\Omega_1(p) - {1\over 2}\Omega_2(p) + {1\over 2}\Omega_3(p)\bigg)
+ O\big(N^{1-\eta}\big),
\leqno(4.4)$$
where 
$$\eqalign{
& \Omega_1(p)
:= S\big({\cal A}_p; {\cal P}(pN), \underline{p}^{1/s'}\big),
\cr\noalign{\bigskip}
& \Omega_2(p)
:= \sum_{\scriptstyle\underline{p}^{1/s'}
\le p_1<\underline{p}^{1/s}
\atop\scriptstyle (p_1, N)=1}
S\big({\cal A}_{pp_1}; {\cal P}(pN), \underline{p}^{1/s'}\big),
\cr
& \Omega_3(p)
:= \mathop{\sum \, \sum \, \sum}_{\scriptstyle
\underline{p}^{1/s'}\le p_1<p_2<p_3<\underline{p}^{1/s}
\atop\scriptstyle (p_1p_2p_3, N)=1}
S\big({\cal A}_{pp_1p_2p_3}; {\cal P}(pp_1N), p_2\big).
\cr}$$
Similar to (5.1), (5.2) and (5.9) of [13], 
we can prove, uniformly for $N\ge 10$ and for $1\le j\le J$, 
$$\sum_p \pi_{[\alpha_{j-1}, \alpha_j)}(p) \Omega_i(p)
\le \big\{\widetilde{\Omega}_i(s, s')+o(1)\big\} 
\Theta\big(N, \pi_{[\alpha_{j-1}, \alpha_j)}\big)
\quad(i=1, 2, 3),$$
where

$$\eqalign{
\widetilde{\Omega}_2(s, s')
& := A(s'),
\cr\noalign{\vskip 1mm}
\widetilde{\Omega}_i(s, s')
& := \int_{1-1/s}^{1-1/s'} {a(s't)\over t(1-t)} \d t,
\cr
\widetilde{\Omega}_3(s, s')
& := 2\max_{\phi\ge 2} 
\mathop{\int\int\int}_{1/s'\le t\le u\le v\le 1/s} 
\omega\bigg({\phi - t - u - v\over u}\bigg)
{\d t \d u \d v\over t u^2 v}.
\cr}$$
Inserting these into (4.4) and 
noticing that
$$
A(s') = 1 + \int_2^{s'-1} {\log(t-1)\over t} \d t,
\qquad
a(s't) = \log(s't-1),
$$
we find that
$$\eqalign{S
& \le  \{1-\Psi_1(s)+o(1)\}
\sum_{1\le j\le J}
\Theta\big(N, \pi_{[\alpha_{j-1}, \alpha_j)}\big)
+ O\big(N^{1-\eta}\big),
\cr
& \le 
\bigg\{
8\big(1-\Psi_1(s)\big)\int_{1/2-s\kappa}^\phi {\d t\over t(1-2t)}
+o(1)\bigg\}
\Theta(N),
\cr}$$
which is equivalent to the required result for the case of $\Psi_1(s)$, 
since 
$$\eqalign{\int_{(1/2-\phi)/\kappa}^s {A(t)\over t(1-2\kappa t)}\d t
& = \int_{(1/2-\phi)/\kappa}^s {\d t\over t(1-2\kappa t)}
\cr
& = \int_{1/2-s\kappa}^\phi {\d t\over t(1-2t)}.
\cr}$$

The case of $\Psi_2(s)$ can be treated in the same way.
The main difference is to use Lemma 4.2 of [13] 
in place of Lemma 4.1 of [13].
We omit the details.
\hfill
$\square$

\goodbreak

\vskip 5mm

\noindent{\bf \S\ 5. Proof of Theorem}

\medskip

Take
$$\kappa_1=1/13.27
\qquad{\rm and}\qquad
\kappa_2=1/8.24,
\leqno(5.1)$$
which satisfy the hypothesis of Lemma 2.2.
Next we  estimate all the terms $\Upsilon_i$ in (2.2).

\smallskip

$1^\circ$
{\it Lower bounds of $\Upsilon_1$ and $\Upsilon_2$}

\smallskip

Write $N^{\kappa}=\underline{1}^{\kappa'}$ 
with $\kappa':=\kappa/(1/2-\delta)$.
By using (4.2) with $\sigma:={\bf 1}_{\{1\}}$
(the charateristic function of $\{1\}$), it follows that
$$\eqalign{\Upsilon_i
& =\Phi(N, {\bf 1}_{\{1\}}, 1/\kappa_i')
\cr
& \ge \big\{a(1/\kappa_i')+h_{k, N_0}(1/\kappa_i')\big\}
\Theta\big(N, {\bf 1}_{\{1\}}\big)
\cr
& \ge \{F_i + o(1)\} \, \Theta(N)
\cr}
\leqno(5.2)$$
with
$$F_i := 8a(1/(2\kappa_i))+8h(1/(2\kappa_i))
\quad(i=1, 2).$$
Write
$$G_2:=8a(1/(2\kappa_i))+8h(1/(2\kappa_2)).$$

$2^\circ$
{\it Upper bounds of $\Upsilon_3$ and $\Upsilon_4$}

\smallskip

We divide the sum $\Upsilon_3$ (resp. $\Upsilon_4$) into subsums
according to
$$\eqalign{
& \hbox{
(a)
$N^{\kappa_1}\le p<N^{1/4}$},
\cr
& \hbox{
(b)
$N^{1/4}\le p<N^{1/2-s_9\kappa_1}$},
\cr
& \hbox{
(c)
$N^{1/2-s_j\kappa_1}\le p<N^{1/2-s_{j-1}\kappa_1}
\;(9\ge j\ge 4)$},
\cr
& \hbox{
(d)
$N^{1/2-s_3\kappa_1}\le p<N^{1/3}$}
\cr}$$
(resp. $N^{\kappa_1}\le p<N^{1/4}$
or $N^{1/4}\le p<N^{1/2-3\kappa_1}$),
where $s_i$ is defined by (3.9).
The contribution of (a) is estimated by Proposition 4.1
and we evaluate (b) (resp. $N^{1/4}\le p<N^{1/2-3\kappa_1}$)
by the classic linear sieve.
The remaining subsums are treated by Proposition 4.4.
It is worth to point out that 
the case (b) requires another kind of treatment
because $\Psi_1(s_{10})=0$ (see Table 3 below).
Thus we obtain
$$\Upsilon_i
\le \{F_i + o(1)\} \, \Theta(N)
\quad(i=3, 4),
\leqno(5.3)$$
where
$$\eqalign{
F_3
& := 8\int_{1/(6\kappa_1)}^{1/(2\kappa_1)-1} 
{A(t)\over t(1-2\kappa_1 t)}\d t
-G_3,
\cr\noalign{\vskip 1mm}
F_4
& := 8\int_3^{1/(2\kappa_1)-1} 
{A(t)\over t(1-2\kappa_1 t)}\d t
-G_4,
\cr}$$
and
$$\eqalign{
G_4
& := 8\int_{(1/(4\kappa_1)}^{1/(2\kappa_1)-1} 
{H(t)\over t(1-2\kappa_1 t)}\d t,
\cr
G_3
& := 8\int_{(1/(4\kappa_1)}^{1/(2\kappa_1)-1} 
{H(t)\over t(1-2\kappa_1 t)}\d t
+ 8\int_{1/(6\kappa_1)}^{s_3} {\Psi_2(s_3)\over t(1-2\kappa t)}\d t
\cr
& \quad
+ 8\sum_{4\le i\le 5}
\int_{s_{i-1}}^{s_i} {\Psi_2(s_i)\over t(1-2\kappa t)}\d t
+ 8\sum_{6\le i\le 9}
\int_{s_{i-1}}^{s_i} {\Psi_1(s_i)\over t(1-2\kappa t)}\d t
\cr}$$

\smallskip

$3^\circ$
{\it Lower bounds of $\Upsilon_5$ and $\Upsilon_6$}

\smallskip

Since $\kappa_1+2\kappa_2=0.318\dots<1/2$, 
Proposition 4.2 yields
$$\Upsilon_5
\ge \{F_5 + o(1)\} \, \Theta(N),
\leqno(5.4)$$
where
$$\eqalign{
F_5 
& := 8
\int_{\kappa_1}^{\kappa_2}
\int_{(1/2-\kappa_2-t)/\kappa_1}^{(1/2-2t)/\kappa_1}
{a(u)\d t\d u\over tu(1-2t-2\kappa_1 u)}
+ G_5
\cr}$$
and
$$G_5
:= 8\int_{\kappa_1}^{\kappa_2} 
\int_{(1/2-\kappa_2-t)/\kappa_1}^{(1/2-2t)/\kappa_1}
{h(u)\d t\d u\over tu(1-2t-2\kappa_1u)}. 
$$

We divide the double sum $\Upsilon_6$ into three subsums
according to
$$\eqalign{
& \hbox{
(a)
$N^{\kappa_1}\le p_1<N^{\kappa_2}\le p_2<N^{1/2-2\kappa_2}$},
\cr
& \hbox{
(b)
$N^{\kappa_1}\le p_1<N^{3\kappa_1/2}$
and 
$N^{1/2-2\kappa_2}\le p_2<N^{1/2-3\kappa_1}$},
\cr
& \hbox{
(c)
$N^{3\kappa_1/2}\le p_1<N^{\kappa_2}$
and 
$N^{1/2-2\kappa_2}\le p_2<N^{1/2-3\kappa_1}$}.
\cr}$$
The first two subsums can be estimated by Proposition 4.3
and the last one by the classic linear sieve.
Thus we obtain
$$\Upsilon_6
\ge \{F_6 + o(1)\} \, \Theta(N),
\leqno(5.5)$$
where
$$F_6 
:= 8
\int_{\kappa_1}^{\kappa_2}
\int_{(3\kappa_1-t)/\kappa_1}^{(1/2-\kappa_2-t)/\kappa_1}
{a(u)\d t\d u\over tu(1-2t-2\kappa_1 u)} 
+ G_6$$
and
$$G_6
:=8\int_{\kappa_1}^{\kappa_2} 
\int_{(2\kappa_2-t)/\kappa_1}^{(1/2-\kappa_2-t)/\kappa_1}
{h(u)\d t\d u\over tu(1-2t-2\kappa_1u)}
+ 8\int_{\kappa_1}^{3\kappa_1/2} 
\int_{(3\kappa_1-t)/\kappa_1}^{(2\kappa_2-t)/\kappa_1}
{h(u)\d t\d u\over tu(1-2t-2\kappa_1u)}. 
$$

\smallskip

$4^\circ$
{\it Upper bounds of $\Upsilon_i$ for $i = 7, 8, 9, 10, 11$}

\smallskip

Clearly the terms 
$\Upsilon_7$, 
$\Upsilon_8$, 
$\Upsilon_9$, 
$\Upsilon_{10}$
and $\Upsilon_{11}$
here
are those terms 
$\Upsilon_7$ (with $\sigma_1=1/2-3\kappa_1$), 
$\Upsilon_9$ (with $\sigma_1=1/3$), 
$\Upsilon_{10}$ (with $\sigma_2=1/2-3\kappa_1$), 
$\Upsilon_{13}$
and $\Upsilon_{14}$
of (9.4) in [13].
Thus (10.10), (10.11), (10.12) of [13] give us the estimates
$$\Upsilon_i
\le \{F_i + o(1)\} \, \Theta(N)
\qquad(i = 7, 8, 9, 10, 11),
\leqno(5.6)$$
where
$$\eqalign{
F_7 
& := 8 \int_2^{2/(1-6\kappa_1)-1} {\log(t-1)\over t} \d t,
\cr
F_8
& := {36\over 5} \int_{\kappa_1}^{1/10} 
{\log(2-3t)\over t (1-t)^2} \d t
+ 8 \int_{1/10}^{1/3} {\log(2-3t)\over t (1-t)} \d t,
\cr
F_9
& := 8
\int_{\kappa_2}^{1/2-3\kappa_1}
{\log\{(1+6\kappa_1-2t)/(1-6\kappa_1)\}\over t (1-t)} \d t,
\cr
F_{10}
& := {36\over 5}
\int_{\kappa_1}^{1/10} {\d t_1\over t_1(1-t_1)}
\int_{t_1}^{\kappa_2} {\d t_2\over t_2^2}
\int_{t_2}^{\kappa_2} {\d t_3\over t_3}
\int_{t_3}^{\kappa_2} \omega\bigg({1-t_1-t_2-t_3-t_4\over t_2}\bigg) 
{\d t_4\over t_4}
\cr
& \quad
+ 8 \int_{1/10}^{\kappa_2} {\d t_1\over t_1}
\int_{t_1}^{\kappa_2} {\d t_2\over t_2^2}
\int_{t_2}^{\kappa_2} {\d t_3\over t_3}
\int_{t_3}^{\kappa_2} \omega\bigg({1-t_1-t_2-t_3-t_4\over t_2}\bigg) 
{\d t_4\over t_4},
\cr
F_{11}
& := {36\over 5}
\int_{\kappa_1}^{1/10} {\d t_1\over t_1(1-t_1)}
\int_{t_1}^{\kappa_2} {\d t_2\over t_2^2}
\int_{t_2}^{\kappa_2} {\d t_3\over t_3}
\int_{\kappa_2}^{1/2-2\kappa_1-t_3} 
\omega\bigg({1-t_1-t_2-t_3-t_4\over t_2}\bigg) 
{\d t_4\over t_4}
\cr
& \quad
+ 8 \int_{1/10}^{\kappa_2} {\d t_1\over t_1}
\int_{t_1}^{\kappa_2} {\d t_2\over t_2^2}
\int_{t_2}^{\kappa_2} {\d t_3\over t_3}
\int_{\kappa_2}^{1/2-2\kappa_1-t_3} 
\omega\bigg({1-t_1-t_2-t_3-t_4\over t_2}\bigg) 
{\d t_4\over t_4},
\cr}$$
and $\omega(t)$ is the Buchstab function (see Lemma 2.10 of [13]).

Inserting (5.2)--(5.6) into (2.2),
we get the following inequality
$$D_{1,2}(N)
\ge \{F(\kappa_1, \kappa_2) + o(1)\} \, \Theta(N),$$
where
$$\eqalign{F(\kappa_1, \kappa_2)
& := \textstyle {1\over 4} 
(3F_1
+ F_2
- F_3
- F_4
+ F_5
+ F_6
- 2F_7
- F_8
- F_9
- F_{10}
- F_{11}).
\cr}$$

\smallskip

$5^\circ$
{\it Numeric computation}

\smallskip

From (3.2) and (3.3), we deduce easily that
$$a(s):=\cases{
0 
& if $\;0<s\le 2$,
\cr\noalign{\vskip 3mm}
\log(s-1) 
& if $\;2<s\le 4$, 
\cr\noalign{\vskip 0,5mm}
\displaystyle
\log(s-1)
+ \int_3^{s-1} {\d t\over t}
\int_2^{t-1} {\log(u-1)\over u} \d u
& if $\;4<s\le 6$, 
\cr\noalign{\vskip 0,5mm}
\displaystyle
\log(s-1)
+ \int_3^{s-1} {\d t\over t}
\int_2^{t-1} {\log(u-1)\over u} \d u
\cr\noalign{\vskip 0,5mm}
\displaystyle
\hskip 5mm
+ \int_5^{s-1} {\d t\over t}
\int_4^{t-1} {\d u\over u}
\int_3^{u-1} {\d v\over v} 
\int_2^{v-1} {\log(w-1)\over w} \d w
& if $\;6<s\le 8$, 
\cr}$$
and
$$A(s):=\cases{
1 
& if $\;0<s\le 3$, 
\cr\noalign{\vskip 2mm}
\displaystyle
1
+ \int_2^{s-1}{\log(t-1)\over t}\d t
& if $\;3<s\le 5$, 
\cr\noalign{\vskip 0,5mm}
\displaystyle
1 
+ \int_2^{s-1}{\log(t-1)\over t}\d t
+ \int_4^{s-1} {\d t\over t}
\int_3^{t-1} {\d u\over u} 
\int_2^{u-1} {\log(v-1)\over v} \d v
& if $\;5<s\le 7$. 
\cr}$$

BY using (3.8), we have
$$G_2\ge 
8\bigg(h(s_{22})+\int_{s_{22}-1}^{1/(2\kappa_2)-1}{H(t)\over t}\d t\bigg)
\ge 0.005283.$$

In order to estimate $G_4$,
we use Table 1 and the decreasing property of $H(s)$ to obtain
$$\eqalign{G_4
= 8\int_{1/(4\kappa_1)}^{1/(2\kappa_1)-1} 
{H(t)\over t(1-2\kappa_1 t)}\d t
\ge 8\sum_{14\le i\le 29} g_4^i H(s_i)
\ge 0.008860
\cr}$$
with
$$\eqalign{
g_4^{14}
& :=\log\bigg({2\kappa_1s_{14}\over 1-2\kappa_1s_{14}}\bigg),
\cr
g_4^i
& :=\log\bigg({s_i(1-2\kappa_1s_{i-1})\over s_{i-1}(1-2\kappa_1s_i)}\bigg)
\quad
(15\le i\le 29).
\cr}$$

With a simpler calculation, we get
$$\eqalign{G_3
& = G_4
+ 8\sum_{3\le i\le 5}g_3^i\Psi_2(s_i)
+ 8\sum_{6\le i\le 9}g_3^i\Psi_1(s_i)
\cr
& \ge 0.039890
\cr}$$
with
$$\eqalign{
g_3^3
& :=\log\bigg({4\kappa_1s_3\over1-2\kappa_1s_3}\bigg),
\cr
g_3^i
& :=\log\bigg({s_i(1-2\kappa_1s_{i-1})\over s_{i-1}(1-2\kappa_1s_i)}\bigg)
\quad(4\le i\le 9).
\cr}$$
Here we have used Table 1 of [13] on the lower bounds for 
$\Psi_2(s_i)$ ($3\le i\le 5$) and $\Psi_1(s_i)$ ($6\le i\le 9$):
$$\displaylines{
\vbox{\tabskip = 0pt\offinterlineskip
\halign{
\vrule # & &\hfil$ $ $#$ $ \!\! $ \hfil & \vrule #\cr
\noalign{\hrule}
height 2mm
&&&&&&&&&&&&&&&&
\cr
& \,\,     i              \,\,        &
& \,\,\,   s_i            \,\,\,      &
& \,\,\,\, s_i'           \,\,\,\,    &
& \,\,\,   \kappa_{1,i}   \,\,\,      &
& \,\,\,   \kappa_{2,i}   \,\,\,      &
& \,\,\,   \kappa_{3,i}   \,\,\,      &
& \quad\,\,\, \Psi_1(s_i) \quad\,\,\, &
& \quad\,\,\, \Psi_2(s_i) \quad\,\,\, &
\cr
height 2mm
&&&&&&&&&&&&&&&&
\cr
\noalign{\hrule}
height 2mm
&&&&&&&&&&&&&&&&
\cr
& 3           &
& 2.3         &
& 4.50        &
& 3.54        &
& 2.88        &
& 2.43        &
&             &
& 0.015247971 &
\cr
height 2mm
&&&&&&&&&&&&&&&&
\cr
\noalign{\hrule}
height 2mm
&&&&&&&&&&&&&&&&
\cr
& 4           &
& 2.4         &
& 4.46        &
& 3.57        &
& 2.87        &
& 2.40        &
&             &
& 0.013898757 &
\cr
height 2mm
&&&&&&&&&&&&&&&&
\cr
\noalign{\hrule}
height 2mm
&&&&&&&&&&&&&&&&
\cr
& 5           &
& 2.5         &
& 4.12        &
& 3.56        &
& 2.91        &
& 2.50        &
&             &
& 0.011776059 &
\cr
height 2mm
&&&&&&&&&&&&&&&&
\cr
\noalign{\hrule}
height 2mm
&&&&&&&&&&&&&&&&
\cr
& 6           &
& 2.6         &
& 3.58        &
&             &
&             &
&             &
& 0.009405211 &
&             &
\cr
height 2mm
&&&&&&&&&&&&&&&&
\cr
\noalign{\hrule}
height 2mm
&&&&&&&&&&&&&&&&
\cr
& 7           &
& 2.7         &
& 3.47        &
&             &
&             &
&             &
& 0.006558950 &
&             &
\cr
height 2mm
&&&&&&&&&&&&&&&&
\cr
\noalign{\hrule}
height 2mm
&&&&&&&&&&&&&&&&
\cr
& 8           &
& 2.8         &
& 3.34        &
&             &
&             &
&             &
& 0.003536751 &
&             &
\cr
height 2mm
&&&&&&&&&&&&&&&&
\cr
\noalign{\hrule}
height 2mm
&&&&&&&&&&&&&&&&
\cr
& 9           &
& 2.9         &
& 3.19        &
&             &
&             &
&             &
& 0.001056651 &
&             &
\cr
height 2mm
&&&&&&&&&&&&&&&&
\cr
\noalign{\hrule}
height 2mm
&&&&&&&&&&&&&&&&
\cr
& 10          &
& 3.0         &
& 3.00        &
&             &
&             &
&             &
& 0           &
&             &
\cr
height 2mm
&&&&&&&&&&&&&&&&
\cr
\noalign {\hrule}}}
\cr}$$
\vskip -2mm
\centerline{Table 3. Lower bounds for $\Psi_1(s_i)$ and $\Psi_2(s_i)$}

\medskip

Similarly
$$\eqalign{G_5
& = 8\int_{\kappa_1}^{\kappa_2} 
\int_{(1/2-\kappa_2-t)/\kappa_1}^{(1/2-2t)/\kappa_1}
{h(u)\d u\d t\over tu(1-2t-2\kappa_1u)}
\cr
& = 8\int_{(1/2-2\kappa_2)/\kappa_1}^{(1/2-\kappa_1-\kappa_2)/\kappa_1} 
h(u)
\log\bigg({2\kappa_2\over 1-2\kappa_2-2\kappa_1u}\bigg)
{\d u\over u(1-2\kappa_1u)}
\cr
& \quad
+ 8\int_{(1/2-\kappa_1-\kappa_2)/\kappa_1}^ {(1/2-2\kappa_1)/\kappa_1}
h(u)
\log\bigg({1-2\kappa_1-2\kappa_1u\over 2\kappa_1}\bigg)
{\d u\over u(1-2\kappa_1u)}
\cr
& \ge 8\sum_{15\le i\le 27} g_5^i h(s_i)
\cr\noalign{\vskip 0,5mm}
& \ge 0.001359
\cr}$$
with
$$\eqalign{
g_5^{15}
& :=\int_{(1/2-2\kappa_2)/\kappa_1}^{s_{15}} 
\log\bigg({2\kappa_2\over 1-2\kappa_2-2\kappa_1u}\bigg)
{\d u\over u(1-2\kappa_1u)},
\cr
g_5^i
& :=\int_{s_{i-1}}^{s_i} 
\log\bigg({2\kappa_2\over 1-2\kappa_2-2\kappa_1u}\bigg)
{\d u\over u(1-2\kappa_1u)}
\qquad
(16\le i\le 20),
\cr
g_5^{21}
& :=\int_{s_{20}}^{(1/2-\kappa_1-\kappa_2)/\kappa_1} 
\log\bigg({2\kappa_2\over 1-2\kappa_2-2\kappa_1u}\bigg)
{\d u\over u(1-2\kappa_1u)}
\cr
& \quad
+\int_{(1/2-\kappa_1-\kappa_2)/\kappa_1}^{s_{21}} 
{\log(1/(2\kappa_1)-1-u)\over u(1-2\kappa_1u)} \d u,
\cr
g_5^i
& :=\int_{s_{i-1}}^{s_i} 
{\log(1/(2\kappa_1)-1-u)\over u(1-2\kappa_1u)} \d u
\qquad
(22\le i\le 26),
\cr
g_5^{27}
& :=\int_{s_{26}}^{(1/2-2\kappa_1)/\kappa_1}
{\log(1/(2\kappa_1)-1-u)\over u(1-2\kappa_1u)} \d u;
\cr}
$$
and
$$\eqalign{
G_6 
& = 8\int_{\kappa_1}^{\kappa_2}
\int_{(2\kappa_2-t)/\kappa_1}^{(1/2-\kappa_2-t)/\kappa_1}
{h(u)\d u\d t\over tu(1-2t-2\kappa_1u)} 
+ 8\int_{\kappa_1}^{3\kappa_1/2} 
\int_{(3\kappa_1-t)/\kappa_1}^{(2\kappa_2-t)/\kappa_1}
{h(u)\d t\d u\over tu(1-2t-2\kappa_1u)}
\cr
& \ge 8\int_2^{(1/2-2\kappa_2)/\kappa_1} 
\log\bigg({\kappa_2(1-2\kappa_1-2\kappa_1u)
\over \kappa_1(1-2\kappa_2-2\kappa_1u)}\bigg)
{h(u)\d u\over u(1-2\kappa_1u)}
\cr
& \quad
+ 8\int_{(1/2-2\kappa_2)/\kappa_1}^{(1/2-\kappa_1-\kappa_2)/\kappa_1}
\log\bigg({(1-2\kappa_1-2\kappa_1u)(1-2\kappa_2-2\kappa_1u)
\over 4\kappa_1\kappa_2}\bigg) 
{h(u)\d u\over u(1-2\kappa_1u)}
\cr
& \ge 8\sum_{1\le i\le 21} g_6^i h(s_i)
\cr\noalign{\vskip 0,5mm}
& \ge 0.060469
\cr}$$
with
$$\eqalign{
g_6^i
& :=\int_{s_{i-1}}^{s_i}
\log\bigg({\kappa_2(1-2\kappa_1-2\kappa_1u)
\over \kappa_1(1-2\kappa_2-2\kappa_1u)}\bigg)
{\d u\over u(1-2\kappa_1u)} 
\qquad
(1\le i\le 14),
\cr
g_6^{15}
& :=\int_{s_{14}}^{(1/2-2\kappa_2)/\kappa_1} 
\log\bigg({\kappa_2(1-2\kappa_1-2\kappa_1u)
\over \kappa_1(1-2\kappa_2-2\kappa_1u)}\bigg)
{\d u\over u(1-2\kappa_1u)}
\cr
& \quad
+\int_{(1/2-2\kappa_2)/\kappa_1}^ {s_{15}}
\log\bigg({(1-2\kappa_1-2\kappa_1u)(1-2\kappa_2-2\kappa_1u)
\over 4\kappa_1\kappa_2}\bigg) 
{\d u\over u(1-2\kappa_1u)},
\cr
g_6^i
& :=\int_{s_{i-1}}^{s_i} 
\log\bigg({(1-2\kappa_1-2\kappa_1u)(1-2\kappa_2-2\kappa_1u)
\over 4\kappa_1\kappa_2}\bigg) 
{\d u\over u(1-2\kappa_1u)}
\qquad
(16\le i\le 20),
\cr
g_6^{21}
& :=\int_{s_{20}}^{(1/2-\kappa_1-\kappa_2)/\kappa_1}
\log\bigg({(1-2\kappa_1-2\kappa_1u)(1-2\kappa_2-2\kappa_1u)
\over 4\kappa_1\kappa_2}\bigg) 
{\d u\over u(1-2\kappa_1u)}.
\cr}
$$

To simplify the computation of $F_{10}$ and $F_{11}$,
we make use of the fact that 
$\omega(t)\le 0.561522$ for $t\ge 3.4$.

Finally
a numerical computation concludes
$$\eqalign{\textstyle F(\kappa_1, \kappa_2)
& \ge \textstyle {1\over 4}\big\{
   3\times 14.900897
+ (9.103015+0.005283)
\cr
& \quad
- (23.652925-0.039890)
- (19.643510-0.008860)
\cr
& \quad
+ (1.654808+0.001359)
+ (3.819092+0.060469)
\cr
& \quad
- 2\times 0.585179
- 5.279581 
- 5.372410 
- 0.104305
- 0.543858\big\}
\cr
& >0.899.
\cr}$$
This completes the proof of Theorem.
\hfill
$\square$

\goodbreak

\vskip 10mm

\centerline{\bf References}

\bigskip

\item{[1]}
{\author E. Bombieri \& H. Davenport},
Small differences between prime numbers,
{\it Proc. Roy. Soc}. Ser. A {\bf 239} (1966), 1--18.

\item{[2]}
{\author Y.C. Cai},
On Chen's theorem (II),
{\it J. Number Theory}, to appear.

\item{[3]}
{\author Y.C. Cai},
Chen's theorem with small primes (Chinese),  
{\it Acta Math. Sinica (Chin. Ser.)} {\bf 48}  (2005), no. 3, 593--598. 

\item{[4]}
{\author Y.C. Cai \& M.G. Lu},
On Chen's theorem,
in: {\it Analytic number theory} (Beijing/Kyoto, 1999), 99--119, 
Dev. Math. {\bf 6}, Kluwer Acad. Publ., Dordrecht, 2002.

\item{[5]}
{\author J.R. Chen},
On the representation of a large even integer
as the sum of a prime and the product of at most two primes,
{\it Sci. Sinica} {\bf 16} (1973), 157--176.

\item{[6]}
{\author J.R. Chen},
On the representation of a large even integer
as the sum of a prime and the product of at most two primes (II),
{\it Sci. Sinica} {\bf 21} (1978), 421--430.

\item{[7]}
{\author J.R. Chen},
Further improvement on the constant in the proposition `1+2':
On the representation of a large even integer
as the sum of a prime and the product of at most two primes (II) 
(in Chinese),
{\it Sci. Sinica} {\bf 21} (1978), 477--494.

\item{[8]}
{\author J.R. Chen},
On the Goldbach's problem and the sieve methods,
{\it Sci. Sinica} {\bf 21} (1978), 701--739.

\item{[9]}
{\author H. Halberstam \& H.-E. Richert},
{\it Sieve Methods},
Academic Press, London, 1974.

\item{[10]}
{\author G.H. Hardy \& J.E. Littlewood},
Some problems of `partitio numerorum' III :
On the expression of a number as a sum of primes,
{\it Acta Math.} {\bf 44} (1923), 1--70.

\item{[11]}
{\author C.D. Pan \& C.B. Pan},
{\it Goldbach Conjecture},
Science Press, Beijing, China, 1992.

\item{[12]}
{\author J. Wu},
Sur la suite des nombres premiers jumeaux,
{\it Acta Arith.} {\bf 55} (1990), 365--394.

\item{[13]}
{\author J. Wu},
Chen's double sieve, Goldbach's conjecture and the twin prime problem,
{\it Acta Arith.} {\bf 114} (2004), 215--273.

\bigskip

Institut Elie Cartan

UMR 7502 UHP-CNRS-INRIA

Universit\'e Henri Poincar\'e (Nancy 1)

54506 Vand\oe uvre--l\`es--Nancy

FRANCE

e--mail: wujie@iecn.u-nancy.fr

\end

$$\eqalign{G_5(\kappa_1, \kappa_2, \kappa_1)
& = 8\int_{\kappa_1}^{\kappa_2} {\d t\over t}
\int_{(1/2-\kappa_2-t)/\kappa_1}^{(1/2-2t)/\kappa_1}
{h(u)\over u(1/2-t-\kappa_1u)} \d u
\cr
& = 8\int_{(1/2-2\kappa_2)/\kappa_1}^{(1/2-\kappa_1-\kappa_2)/\kappa_1} 
{h(u)\over u}\d u
\int_{1/2-\kappa_2-\kappa_1u}^{(1/2-\kappa_1u)/2}
{\d t\over 1/2-t-\kappa_1u}
\cr
& \quad
+ 8\int_{(1/2-\kappa_1-\kappa_2)/\kappa_1}^ {(1/2-2\kappa_1)/\kappa_1}
{h(u)\over u}\d u
\int_{\kappa_1}^{(1/2-\kappa_1u)/2}
{\d t\over 1/2-t-\kappa_1u}
\cr
& = 8\int_{(1/2-2\kappa_2)/\kappa_1}^{(1/2-\kappa_1-\kappa_2)/\kappa_1} 
{h(u)\over u}\log\bigg({4\kappa_2\over 1-2\kappa_1u}\bigg)\d u
\cr
& \quad
+ 8\int_{(1/2-\kappa_1-\kappa_2)/\kappa_1}^ {(1/2-2\kappa_1)/\kappa_1}
{h(u)\over u}
\log\bigg({2-4\kappa_1-4\kappa_1u\over 1-2\kappa_1u}\bigg)\d u
\cr}$$

$$\eqalign{G_6(\kappa_1, \kappa_2, 1/2-3\kappa_1, \kappa_1)
& =8\int_{\kappa_1}^{\kappa_2} {\d t\over t}
\int_{(3\kappa_1-t)/\kappa_1}^{(1/2-\kappa_2-t)/\kappa_1}
{h(u)\over u(1/2-t-\kappa_1u)} \d u
\cr
& =8\int_{3-\kappa_2/\kappa_1}^{2} {h(u)\over u} \d u
\int_{\kappa_1(3-u)}^{\kappa_2}{\d t\over 1/2-t-\kappa_1u}
\cr
& \quad
+ 8\int_{2}^{(1/2-2\kappa_2)/\kappa_1} {h(u)\over u} \d u
\int_{\kappa_1}^{\kappa_2}{\d t\over 1/2-t-\kappa_1u}
\cr
& \quad
+ 8\int_{(1/2-2\kappa_2)/\kappa_1}^{(1/2-\kappa_1-\kappa_2)/\kappa_1}
 {h(u)\over u} \d u
\int_{\kappa_1}^{1/2-\kappa_2-\kappa_1u}{\d t\over 1/2-t-\kappa_1u}
\cr
& =8\int_{3-\kappa_2/\kappa_1}^{2} 
{h(u)\over u}
\log\bigg({1-6\kappa_1\over 1-2\kappa_2-2\kappa_1u}\bigg) \d u
\cr
& \quad
+ 8\int_{2}^{(1/2-2\kappa_2)/\kappa_1} 
{h(u)\over u} 
\log\bigg({1-2\kappa_1-2\kappa_1u\over 1-2\kappa_2-2\kappa_1u}\bigg)\d u
\cr
& \quad
+ 8\int_{(1/2-2\kappa_2)/\kappa_1}^{(1/2-\kappa_1-\kappa_2)/\kappa_1}
{h(u)\over u}
\log\bigg({1-2\kappa_1-2\kappa_1u\over 2\kappa_2}\bigg) \d u
\cr}$$

\item{[2]}
{\author E. Bombieri, J.B. Friedlander \& H. Iwaniec},
Primes in arithmetic progressions to large moduli,
{\it Acta Math}. {\bf 156} (1986), 203--251.

\bye